# Transformation, Identification, and Inversion of *Goldberg-Coxeter* Fullerenes

LI, Shaoqing

Anhui Jianzhu university, Hefei, Anhui, China.

**Abstract:** It is difficult to identify a *G-C* fullerene directly from its dimensions as its lattice is not proportional to that of its archetype in general, although they have the same three-dimensional shape. In this paper, the area scale factor of a *G-C* fullerene is proved to be an integer, which can be calculated from its dimensions. All the *G-C* transformations are *k*-inflations that can be easily identified and inversed, primary transformations whose area scale factors are prime numbers, or composite transformations whose area scale factors are the product of those of its sub-transformations. As the result, a method to identify any *G-C* fullerenes according to the area scale factor was presented.

## 1. Introduction

The *Goldberg-Coxeter* transformation is a global transformation for fullerene graphs. Given any fullerene, it is possible to construct an infinite series of larger ones with the same three-dimensional shape through the *G-C* transformations [1, 2]. The constructed fullerenes are called *Goldberg-Coxeter* fullerenes. The classic *G-C* construction for icosahedral fullerenes, the leapfrog transformation, and the *k*-inflation are all special cases of the *G-C* transformation [1].

Although a *G-C* fullerene has the same shape as its archetype, its lattice is not proportional to that of its archetype in general. Except for those special cases, it is difficult to identify a *G-C* fullerene directly from its dimensions.

With fullerene unfolded onto the Eisenstein plane, a *G-C* transformation can be expressed as multiplication in Euclid's algorithm [1]. Similarly, a *G-C* transformation can be inverted by dividing instead of multiplying, and we can even use Euclid's algorithm to find out whether a particular fullerene is a *G-C* fullerene and to find its archetype [1]. To do such inversion and identification, we must get reasonable candidates for the *G-C* transformation first.

The shape of a fullerene can be characterized by its master polyhedron [3-5]. With a *G-C* transformation, the master polyhedron is scaled up. The scale factor of the face areas is the same as that of the number of carbon atoms, which is an integer. As the face areas and the edge squares of the master polyhedron are all integers. the area scale factor is a common factor of those integers. With the area scale factor, we can get the candidates for the *G-C* transformation.

For the *G-C* fullerene without *k*-inflation, each prime factor of the area scale factor corresponds to a possible sub-transformation. Therefore, we can identify the *G-C* fullerene without k-inflation by check each prime factor in the area scale factor by inversion.

In this paper, the faces of the master polyhedron are called the facets of the fullerene, and the edges of the master polyhedron are called the master edges of the fullerene. The remainder of the paper is as follows. In section 2, the expression of *G-C* transformation in Euclid's algorithm was present, and the identification and inversion method of *k*-inflation fullerene and leapfrog fullerene were presented. In section 3, the area scale factor of a *G-C* transformation was proved to be the common factor of the facet areas and master edge squares. In section 4, the transformation candidates according to the area scale factor was studied. In section 5, the method of identification and inversion of any *G-C* fullerenes was presented.

## 2. Transformation and simple identifications

A *G-C* transformation is to replace each triangular face of the dual of a fullerene (i.e. the graphene patch (1,0) of the fullerene as Figure 1) with another graphene patch (*k, l*), constructing a larger one with the same three-dimensional shape. The transformation was first used to construct icosahedral fullerenes whose archetype is C20.

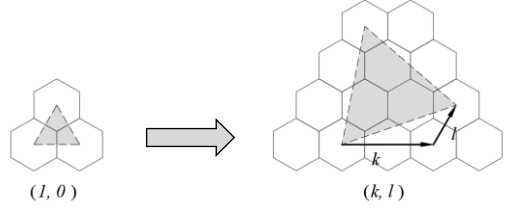

**Fig. 1**. The principle of the *G-C* transformation

Fullerene can be unfolded onto a hexagonal latticed plane which is congruent to an Eisenstein plane characterized by a triangular lattice. The Eisenstein plane is a complex plan with unit vectors in the six directions denoted with 1, $\omega$, $\omega^2$, $\omega^3$, $\omega^4$, $\omega^5$, where $\omega = e^{i2\pi/6}$ is the complex root of equation $\omega^6 = 1$. As $\omega^3 = -1$, $\omega^2 = \omega - 1$, each vector in the Eisenstein plane can be express as ($a+b\omega$), with (*a, b*) known as Eisenstein integers.

In the Eisenstein plane, a *G-C* transformation (*k, l*) can be expressed as

$$G' = (k + l\omega) \cdot G$$

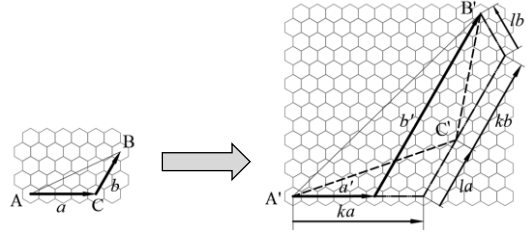

**Fig. 2**. Dimensional change in *G-C* transformation

In Figure 2,
$$a'+b'\omega = (k+l\omega) \cdot (a+b\omega)$$
$$= (ak - bl) + (al + bk + bl)\omega$$

Therefore,

$$\begin{pmatrix} a' \\ b' \end{pmatrix} = \begin{pmatrix} k & -l \\ l & k+l \end{pmatrix} \begin{pmatrix} a \\ b \end{pmatrix} \quad (1)$$

In formula (1), (*a, b*), and (*a', b'*) are all Eisenstein integers. If a pair of *Eisenstein integers* are all great than or equal to zero, they are equal to the Coxeter coordinates. Setting initial Eisenstein integers as the Coxeter coordinates, formula (1) can be used for the calculation of the transformation of the master edges one by one without unfolding the fullerene altogether. The results can be easily transformed to the Coxeter coordinates if needed.

The formula (1) can also be used for the inversion and identification of *G-C* fullerenes. If a pair of (*k, l*) can be found so that the inversion results are all integers, the fullerene is a *G-C* fullerene, and the inversion results correspond to the dimensions of its archetype.

It is not necessary to calculate Coxeter coordinates of all the master edges for transformation,

identification, or inversion. A fullerene has 20 dimensions [3]. Coxeter coordinates for 10 independent master edges are enough for a fullerene.

**Theorem 1.** The fullerene is a *G-C* fullerene if the Coxeter coordinates of its master edges have a common factor great than 1.

**Proof.** According to equation (1), if $l = 0$, then

$$\begin{cases} a' = ka \\ b' = kb \end{cases} \quad (2)$$

This means the lattice is scaled up without rotation in the transformation. Such *Goldberg-Coxeter* transformation is called *k*-inflation.

If the Coxeter coordinates of its master edges have a common divisor $k$ great than 1, the fullerene can result from *k*-inflation. With the factor $k$, the archetype can easily be got according to formula (2). Therefore, the fullerene is a *G-C* fullerene of inflation.

□

**Theorem 2.** The fullerene is a leapfrog fullerene if and only if each pair of Coxeter coordinates of the master edges are congruent modulo 3.

**Proof.** The *G-C* transformation (1, 1) is also called the leapfrog transformation. According to Formula (1),

$$\begin{cases} a = a' + (b' - a')/3 \\ b = (b' - a')/3 \end{cases} \quad (3)$$

So, if and only if two Coxeter coordinates of any edge are congruent modulo 3, $(a, b)$ have integer values, and the fullerene is a leapfrog fullerene.

□

## 3. Area scale factor concealed in the dimensions

The equivalent area (area for short) of a graphene patch is equal to the corresponding number of carbon atoms. This means the area of the graphene patch (1, 0) is 1. With **i**, **j** as the unit vector in the horizon and 60° respectively, the area of a triangular graphene patch expanded by vectors $(a + b\omega)$ and $(c + d\omega)$ in the Eisenstein plane can be expressed as.

$$N_{abcd} = \frac{(a\mathbf{i} + b\mathbf{j}) \times (c\mathbf{i} + d\mathbf{j})}{\mathbf{i} \times \mathbf{j}} = \begin{vmatrix} a & b \\ c & d \end{vmatrix} \quad (4)$$

The formula is just the same as that for the area of a parallelogram in the Cartesian plane.

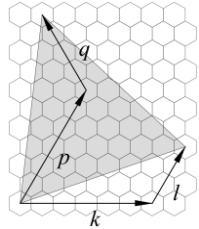

**Fig3**. Coxeter coordinates for a facet.

To a fullerene facet characterized with Coxeter coordinates $(k, l)$ and $(p, q)$ as in Figure 3, the vectors corresponding to the two sides of the facet in the Cartesian plane are $(k+l\omega)$ and $(p\omega+q\omega^2)$ respectively. With $\omega^2 = \omega - 1$, then

$$p\omega + q\omega^2 = p\omega + q(\omega - 1)$$
$$= -q + (p+q)\omega$$

So, according to formula (4), the area of the facet can be calculated as

$$N = \begin{vmatrix} k & l \\ -q & p+q \end{vmatrix} = kp + kq + lq \qquad (5)$$

For the graphene patch *(k, l)* in figure 1, $p = k$ and $q = l$. So its area is

$$N_{kl} = \begin{vmatrix} k & l \\ -l & k+l \end{vmatrix} = k^2 + kl + l^2 \qquad (6)$$

According to formula (6), the area of a triangular patch is equal to the square of its edge. To the transformation *(k, l)*, $N_{kl}$ is the area scale factor.

**Theorem 3.** For a *G-C* fullerene, the areas of its facets and the squares of its master edges are all integers with a common factor as the area scale factor of the *G-C* transformation.

**Proof**. According to formula (5-6), the facet areas and the edge squares of fullerene are all integers. With the *G-C* transformation, they have a common scale factor which is also an integer according to formula (6). Therefore, the area scale factor is a common factor of the facet areas and edge squares of the *G-C* fullerene.

□

## 4. Candidates for *G-C* transformation

With the area scale factor $N_{kl}$, we can calculate the candidates for the *G-C* transformation according to formula (6). The *G-C* transformations with area scale factor $N_{kl} < 100$ are listed in Table 1. For simplicity, the transformation *(k, l)* in the table also represents the chiral one if present, which has the same area scale factor.

**Table 1.** The correspondence between *G-C* transformations *(k, l)* and area scale factors $N_{kl} < 100$.

| *(k, l)* | $N_{kl}$ | Notation | *(k, l)* | $N_{kl}$ | Notation | *(k, l)* | $N_{kl}$ | Notation |
|---|---|---|---|---|---|---|---|---|
| (1, 1) | 3 | *Primary* | (5, 1) | 31 | *Primary* | (8, 0) | 64 | *Inflation* |
| (2, 0) | 4 | *Inflation* | (6, 0) | 36 | *Inflation* | (7, 2) | 67 | *Primary* |
| (2, 1) | 7 | *Primary* | (4, 3) | 37 | *Primary* | (8, 1) | 73 | *Primary* |
| (3, 0) | 9 | *Inflation* | (5, 2) | 39 | *(1, 1) · (3, 1)* | (5, 5) | 75 | *(1, 1) · (5, 0)* |
| (2, 2) | 12 | *(1, 1) · (2, 0)* | (6, 1) | 43 | *Primary* | (6, 4) | 76 | *(3, 2) · (2, 0)* |
| (3, 1) | 13 | *Primary* | (4, 4) | 48 | *(1, 1) · (4, 0)* | (7, 3) | 79 | *Primary* |
| (4, 0) | 16 | *Inflation* | (5, 3) | 49 | *(2, 1) · (2, 1)* | (9, 0) | 81 | *Inflation* |
| (3, 2) | 19 | *Primary* | (7, 0) | 49 | *Inflation* | (8, 2) | 84 | *(4, 1) · (2, 0)* |
| (4, 1) | 21 | *(1, 1) · (2, 1)* | (6, 2) | 52 | *(3, 1) · (2, 0)* | (6, 5) | 91 | *(2, 1) · (3, 1)* |
| (5, 0) | 25 | *Inflation* | (7, 1) | 57 | *(1, 1) · (2, 3)* | (9, 1) | 91 | *(2, 1) · (3, 1)* |
| (3, 3) | 27 | *(1, 1) · (3, 0)* | (5, 4) | 61 | *Primary* | (7, 4) | 93 | *(1, 1) · (5, 1)* |
| (4, 2) | 28 | *(2, 0) · (2, 1)* | (6, 3) | 63 | *(2, 1) · (3, 0)* | (8, 3) | 97 | *Primary* |

As shown in table 1, Many *G-C* transformation are composite transformations which can be decomposed into consecutive transformations. The composition of two consecutive transformations can be expressed with a multiplication of the Eisenstein integers as follows. The composition results in table 1 have been converted into Coxeter coordinates.

$$(a,b)\cdot(c,d) = ac + (ad+bc)\omega + bd\omega^2$$
$$= (ac-bd,\ bc+ad+bd) \quad (7)$$

According to formula (7), the composite transformation is insensitive to the sequence of the sub transformations but sensitive to their chirality. As the result, there are four *G-C* transformations (two pairs of chiral *G-C* transformations) corresponding to the area scale factor 91, and three *G-C* transformations (a pair of chiral *G-C* transformations and a *k*-inflation) corresponding to the area scale factor 49.

A composite *G-C* transformation is determined by its sub transformations. There are mainly two kinds of basic sub transformations: *k*-inflations and primary *G-C* transformations. If a *G-C* transformation is not a *k*-inflation and cannot be decomposed, it is called a primary *G-C* transformation.

The area scale factors of the primary *G-C* transformations in table 1 are all prime numbers. We have tested that area scale factors in 100-200 for the primary *G-C* transformations are also prime numbers. Then, we have the following conjecture.

**Conjecture 1**: The area scale factors of the primary *G-C* transformations are all prime numbers.

The area scale factor of a composite *G-C* transformation is equal to the product of those of its sub-transformations. With conjecture 1 we can get the following conjecture.

**Conjecture 2**: To the *G-C* fullerene without *k*-inflation, each prime factor of the area scale factor corresponds to a primary sub-transformation.

The prime factor 3 corresponds to the leapfrog transformation. Except for the prime factor 3, each prime factor corresponds to a pair of chiral transformations.

## 5. Identification of *G-C* fullerenes

A *G-C* fullerene can be identified by searching for all possible basic transformations.

First, calculate the greatest common factor of its dimensions. If the factor is great than 1, the fullerene is a *k*-inflation fullerene. Dividing all dimensions with the common factor, we can get a fullerene without *k*-inflation.

Second, for a fullerene without *k*-inflation, calculate its area scale factor as section 3. If the area scale factor is great than 1, each prime factor of the area scale factor corresponds to a possible sub transformation.

Third, test the candidates for each sub transformation by inversion, and identify all the reasonable sub transformations.

With this method, we can identify any *G-C* fullerene, and get its transformation (including all sub-transformations) and the ultimate archetype.

Although this identification method is based on conjecture 1 and conjecture 2, there are no problem if the prime factors of area scale factor are all less than 200. It is enough in practice.